\documentclass[11pt]{amsart}
\usepackage{amsmath,amssymb}
\usepackage[all,cmtip]{xy}

\begin{document}

\newtheorem{thm}{Theorem}[section]
\newtheorem{lem}[thm]{Lemma}
\newtheorem{prop}[thm]{Proposition}
\newtheorem{cor}[thm]{Corollary}
\newtheorem{defn}[thm]{Definition}
\newtheorem*{remark}{Remark}

\numberwithin{equation}{section}

\newcommand{\Z}{{\mathbb Z}} %cph changed from \mathbf
\newcommand{\Q}{{\mathbb Q}}
\newcommand{\R}{{\mathbb R}}
\newcommand{\C}{{\mathbb C}}
\newcommand{\N}{{\mathbb N}}
\newcommand{\FF}{{\mathbb F}}
\newcommand{\fq}{\mathbb{F}_q}
\newcommand{\feq}{\overline{\mathbb F}_q}

\newcommand{\rmk}[1]{\footnote{{\bf Comment:} #1}}

\renewcommand{\mod}{\;\operatorname{mod}}
\newcommand{\ord}{\operatorname{ord}}
\newcommand{\TT}{\mathbb{T}}
\renewcommand{\i}{{\mathrm{i}}}
\renewcommand{\d}{{\mathrm{d}}}
\renewcommand{\^}{\widehat}
\newcommand{\HH}{\mathbb H}
\newcommand{\Vol}{\operatorname{vol}}
\newcommand{\area}{\operatorname{area}}
\newcommand{\tr}{\operatorname{tr}}
\newcommand{\norm}{\mathcal N} % norm =(\frac{ n+\sqrt{n^2-4}} 2)^2
\newcommand{\intinf}{\int_{-\infty}^\infty}
\newcommand{\ave}[1]{\left\langle#1\right\rangle} %  average
\newcommand{\Var}{\operatorname{Var}}
\newcommand{\Prob}{\operatorname{Prob}}
\newcommand{\sym}{\operatorname{Sym}}
\newcommand{\disc}{\operatorname{disc}}
\newcommand{\CA}{{\mathcal C}_A}
\newcommand{\cond}{\operatorname{cond}} % conductor
\newcommand{\lcm}{\operatorname{lcm}}
\newcommand{\Kl}{\operatorname{Kl}} %Kloosterman sum
\newcommand{\leg}[2]{\left( \frac{#1}{#2} \right)}  % Legendre symbol

\newcommand{\sumstar}{\sideset \and^{*} \to \sum}

\newcommand{\LL}{\mathcal L} %L-function of u
\newcommand{\sumf}{\sum^\flat}
\newcommand{\Hgev}{\mathcal H_{2g+2,q}}
\newcommand{\USp}{\operatorname{USp}}
\newcommand{\conv}{*}
\newcommand{\dist} {\operatorname{dist}}
\newcommand{\CF}{c_0} % Fejer constant
\newcommand{\kerp}{\mathcal K}

\newcommand{\fs}{\mathfrak S}
\newcommand{\rest}{\operatorname{Res}} % resultant
\newcommand{\af}{\mathbb A} % affine line

\title[Chowla's conjecture for the rational function field]
{The autocorrelation of the M\"obius function and Chowla's
conjecture for the rational function field}
\author{Dan Carmon  and Ze\'ev Rudnick}
\address{Raymond and Beverly Sackler School of Mathematical Sciences,
Tel Aviv University, Tel Aviv 69978, Israel}

\date{\today}
\thanks{The research leading to these results has received funding from the European
Research Council under the European Union's Seventh Framework Programme
(FP7/2007-2013) / ERC grant agreement n$^{\text{o}}$ 320755.}
%\thanks{This work was   supported by  the Israel Science Foundation (grant No.
%1083/10).}
\begin{abstract}
We prove a function field version of Chowla's conjecture on the
autocorrelation of the M\"obius function in the limit of a large
finite field.
\end{abstract}
\maketitle

\section{Introduction}

 There is a well-known equivalence between the Riemann Hypothesis (RH) %on the zeros of the Riemann zeta function
 and square-root cancellation in sums of the M\"obius function $\mu(n)$,
 namely, RH is equivalent to $\sum_{n\leq N} \mu(n)=O(N^{1/2+o(1) })$.
 This sum measures the correlation between $\mu(n)$ and the constant function.
 Recent studies have explored the correlation between $\mu(n)$ and other sequences, see \cite{Green Tao, CS, BSZ}.
 %For instance, Green and Tao \cite{Green Tao} show that the M\"obius function does not correlate with any "nil-sequence".
 Sarnak \cite{Sarnak} showed that $\mu(n)$  does not correlate
 with  any ``deterministic'' (i.~e., zero entropy) sequence, assuming an old conjecture of
 Chowla \cite{Chowla} on the auto-correlation of the M\"obius function, which asserts that
 given an $r$-tuple of distinct integers $\alpha_1,\dots,\alpha_r$ and
 $\epsilon_i\in \{1,2\}$, not all even,
 then  
\begin{equation}\label{chowla's conjecture}
\lim_{N\to \infty} \frac 1N \sum_{n\leq N}\mu(n+\alpha_1
)^{\epsilon_1}\cdot \dots \cdot \mu(n+\alpha_r)^{\epsilon_r}=0
\end{equation}
Note that the number of nonzero summands here, that is the number of
$n\leq N$ for which $n+\alpha_1,\dots n+\alpha_r$ are all
square-free, is asymptotically $c(\alpha)N$, where $c(\alpha)>0$ if
the numbers $\alpha_1,\dots, \alpha_r$ do not contain a complete
system of residues modulo $p^2$ for every prime $p$ \cite{Mirsky},
so that \eqref{chowla's conjecture} is about non-trivial
cancellation in the sum.

Chowla's conjecture \eqref{chowla's conjecture} seems intractable at
this time, the only known case being $r=1$ where it is equivalent
with the Prime Number Theorem. Our goal in this note is to prove a
function field version of Chowla's conjecture.

Let $\fq$ be a finite field  of $q$ elements, and $\fq[x]$ the
polynomial ring over $\fq$. The M\"obius function of a nonzero
polynomial $F\in \fq[x]$ is defined to be $\mu(F)=(-1)^r$ if
$F=cP_1\dots P_r$ with $0\neq c\in \fq$ and $P_1,\dots, P_r$ are
distinct monic irreducible polynomials, and $\mu(F)=0$ otherwise.

Let $M_n\subset \fq[x]$ be the set of monic polynomials of degree
$n$ over $\fq$, which is of size $\#M_n = q^n$. The number of
square-free polynomials in $M_n$ is, for $n>1$, equal to
$q^n-q^{n-1}$ \cite[Chapter 2]{Rosen}. Hence, given $r$  distinct polynomials
$\alpha_1,\dots, \alpha_r \in \fq[x]$, with $\deg \alpha_j<n$, the
number of $F\in M_n$ for which all of $F(x)+\alpha_j(x)$ are
square-free is $q^n + O(rq^{n-1}) $ as $q\to \infty$.

For $r>0$,  distinct polynomials  $\alpha_1,\dots, \alpha_r \in
\fq[x]$, with $\deg \alpha_j<n$ and
 $\epsilon_i\in \{1,2\}$, not all even, set
\begin{equation}
  C(\alpha_1,\dots, \alpha_r;n):=
\sum_{F\in M_n} \mu(F+\alpha_1)^{\epsilon_1}\dots
\mu(F+\alpha_r)^{\epsilon_r}
\end{equation}
For $r=1$ and $n>1$ we have $\sum_{F\in M_n} \mu(F) = 0$ \cite[Chapter 2]{Rosen}.
  For $n=1$ we have $\mu(F)\equiv -1$ and the sum equals $(-1)^{\sum
  \epsilon_j} q^n$.
  %there is no cancellation in the sum.
 For $n>1$, $r>1$ we show:
\begin{thm}\label{main thm}
Fix $r>1$ and assume that $n>1$ and $q$ is odd.
%Then for fixed $n$,
Then for any choice of distinct polynomials  $\alpha_1,\dots,
\alpha_r \in \fq[x]$, with $\max \deg \alpha_j<n$,  and
$\epsilon_i\in \{1,2\}$, not all even
\begin{equation}
| C(\alpha_1,\dots, \alpha_r;n)| \leq 2 rn q^{n-\frac 12} + 3 rn^2
q^{n-1}
\end{equation}
\end{thm}
%\begin{verbatim}
%at this point we need to assume n is at least 5
%\end{verbatim}
%An examination of the argument shows that $c_{r,n} \ll r^2n^2$.
Thus for fixed $n>1$,
\begin{equation}
\lim_{q\to \infty} \frac 1{\#M_n}\sum_{F\in M_n}
\mu(F+\alpha_1)^{\epsilon_1}\dots \mu(F+\alpha_r)^{\epsilon_r} =0
\end{equation}
under the assumption of Theorem~\ref{main thm}, giving an analogue
of Chowla's conjecture \eqref{chowla's conjecture}.

Our starting point is Pellet's formula , see e.g. \cite[Lemma 4.1]{Conrad}, which asserts
that for the polynomial ring $\fq[x]$ with $q$ odd (hence the restriction on the parity of $q$ in Theorem~\ref{main thm}), the M\"obius
function $\mu(F)$ can be computed in terms of the discriminant $\disc(F)$ of
$F(x)$ as
\begin{equation}\label{Pellet}
\mu(F) = (-1)^{\deg F} \chi_2(\disc(F))
\end{equation}
where $\chi_2$ is the quadratic character of $\fq$. That will allow
us to express $C(\alpha_1,\dots,\alpha_r;n)$ as a character sum and
to estimate it.

\section{Reduction to a counting problem}

\subsection{Character sums}

We use Pellet's formula \eqref{Pellet} to write
\begin{equation}\label{transform C}
C(\alpha_1,\dots, \alpha_r;n) = (-1)^{nr} \sum_{F\in M_n}
\chi_2\Big(\disc(F+\alpha_1)^{\epsilon_1}\dots
\disc(F+\alpha_r)^{\epsilon_r} \Big)
\end{equation}

Since $\disc(F)$ is polynomial in the coefficients of $F$,
\eqref{transform C} is an $n$-dimensional character sum; will will
estimate it by trivially bounding all but one variable. We single
out the constant term $t:= F(0)$ of $F\in M_n$ and write $F(x) = f
(x)+t$, with
\begin{equation}
f(x) = x^n+a_{n-1} x^{n-1}+\dots +a_1 x
\end{equation}
and set
\begin{equation}
D_f(t) := \disc(f(x)+t)%,\quad D_a(t+1)  = \disc(f_{a,t+1})
\end{equation}
which is a polynomial of degree $n-1$ in $t$. Therefore we have
\begin{equation}
\left| C(\alpha_1,\dots, \alpha_r;n) \right|  \leq \sum_{a\in \fq^{n-1}}
\left|\sum_{t\in \fq} \chi_2\Big(D_{f+\alpha_1}(t)^{\epsilon_1}\dots
D_{f+\alpha_r}(t)^{\epsilon_r}\Big) \right|
\end{equation}

 We use Weil's theorem (the Riemann Hypothesis for curves over a finite field), which implies that for a
  polynomial $P(t)\in \fq[t]$ of positive degree, which is not proportional to a
square of another polynomial, we have \cite[\S 2]{Schmidt}
\begin{equation}\label{weil}
\left| \sum_{t\in \fq}\chi_2(P(t)) \right| \leq (\deg P-1) q^{1/2},
\quad P(t)\neq cH^2(t)
\end{equation}

For us, the relevant polynomial is $P(t) =
D_{f+\alpha_1}(t)^{\epsilon_1}\dots D_{f+\alpha_r}(t)^{\epsilon_r}$,
which has degree $\leq 2r(n-1)$. Instead of requiring that it not be
proportional to a square, we impose the stronger requirement that for some $i$ with
$\epsilon_i$ odd, $D_{f+\alpha_i}(t)$ has positive degree and is
squarefree and that for all $j$ such that $j\neq i$, $D_{f+\alpha_i}(t)$ and
$D_{f+\alpha_j}(t)$ are coprime. We denote the set of coefficients $a$
satisfying the stronger condition by $G_n$ (the ``good'' $a$'s, where we can apply \eqref{weil}), and let
$G^c_n=\fq^{n-1} \backslash G_n$ be the complement of $G_n$, where we use the trivial bound $q$ on the character sum.
Thus we deduce that we can bound
\begin{equation}
\begin{split}
\left | C(\alpha_1,\dots, \alpha_r;n)  \right| &\leq \sum_{a\in G_n}
(2r(n-1)-1)\sqrt{q} +\sum_{a\notin G_n} q \\
&\leq (2r(n-1)-1)q^{n-\frac 12} + q\#G_n^c \end{split}
\end{equation}
where we have used the trivial bound $\#G_n\leq q^{n-1}$ for the first part of the sum.
 Theorem~\ref{main
thm} will follow from
\begin{prop}\label{Prop Gcomp}
Assume that $n>1$ and $\max \deg \alpha_j<n$. Then
  $$\#G^c_n \leq 3 r n^2 q^{n-2}$$
\end{prop}

\subsection{Bounding $\#G_n^c$}
We can write $G^c_n\subset A_n\cup B_n$ where:
\begin{itemize}
\item
 $A_n=A_{n,i}$ is the set of those $a\in
\fq^{n-1}$ for which $D_{f+\alpha_i}(t)$ is either a constant or is
not square-free, that is
%\footnote{this is a bit confusing because
%$D_f(t) = \disc f_{a,t}$ was itself defined as a discriminant...}
\begin{equation}
A_{n}  = \{a\in \fq^{n-1}:  D_{f+\alpha_i}(t)  \mbox{ is constant or
} \disc(D_{f+\alpha_i}) =0 \}
\end{equation}

\item
$B_n=\cup_{j\neq i} B(j)$ where $B(j)$ are those $a$'s for which
  $D_{f+\alpha_i}(t)$ and $D_{f+\alpha_j}(t)$ have a common zero, which
  can be written as the vanishing
of their resultant
\begin{equation}
B(j) = \{a\in \fq^{n-1}: \rest(D_{f+\alpha_i}(t), D_{f+\alpha_j}(t)
) = 0 \}
\end{equation}
\end{itemize}

What is crucial is that  $A_n$ and each $B(j)$ are the zero sets of
a polynomial equation in the coefficients $a$;
%, the equation being
%independent of $\fq$ ("defined over $\Z$");
this is a key property of the discriminant and the resultant.

We will need the following elementary but useful uniform upper bound
on the number of zeros of polynomials (cf \cite[\S 4, Lemma
3.1]{Schmidt}):
\begin{lem}\label{lem:schmidt}
Let $h(X_1,\dots, X_m) \in \fq[X_1,\dots, X_m]$ be a non-zero
polynomial of total degree at most $d$. Then the number of zeros  of
$h(X_1,\dots, X_m)$ in $\fq^m$ is at most
\begin{equation}
\#\{x\in \fq^m: h(x)=0\} \leq d q^{m-1}
\end{equation}
\end{lem}

 As we will see below (see \S~\ref{sec:disc formulae}), the equation defining $A_n$ has total degree
 $3(n-1)(n-2)$ in the coefficients $a_1,\dots,a_{n-1}$, and the
 equation defining  $B(j)$ has total degree $\leq 3(n-1)^2$.
Therefore, by Lemma~\ref{lem:schmidt},  if we show that the
equations defining $A_n$, $B(j)$   are not identically zero, then we
will have proved
\begin{equation}\label{An small}
\#A_n \leq 3n^2q^{n-2}
\end{equation}
and
\begin{equation}\label{Bn small}
 \#B_n\leq 3(r-1)n^2q^{n-2}
 \end{equation}
This immediately gives Proposition~\ref{Prop Gcomp}.

In order to show that a polynomial $h \in \fq[X_1,\dots, X_m]$ is
not identically zero, we may instead consider it as a polynomial
defined over $\feq$, the algebraic closure of $\fq$. In this context
we can investigate the zero set $Z_h = \{a \in \feq^m : h(a) = 0\}$,
which is a subvariety of the affine space $\af^m$. The polynomial
$h$ is not identically zero if and only $Z \neq \af^m$. This shall
be our main tool in the following sections.

\subsection{ Resultant and discriminant formulas}
\label{sec:disc formulae}

The discriminant $\disc(F)$ of a polynomial
$F(x)=a_nx^n+a_{n-1}x^{n-1}+\dots+a_0$, $a_n\neq 0$, is given in
term of its roots $r_1,\dots, r_n$ in the algebraic closure $\feq$
as $\disc F = a_n^{2n-2}\prod_{i<j}(r_i-r_j)^2$, and is a
homogeneous polynomial with integer coefficients in $a_0,\dots
,a_n$, with degree of homogeneity $2n-2$, has total degree $2n-2$,
and has degree $n-1$ as a polynomial in $a_0$. Moreover if $a_i$ is
regarded as having degree $i$ then $\disc(F)$ is homogeneous of
degree $n(n-1)$, that is for every monomial $c_r \prod_i a_i^{r_i}$ in $\disc(F)$, 
\begin{equation}\label{weighted homog}
\sum_i ir_i=n(n-1)
\end{equation}

The resultant of two polynomials $F(x)=a_nx^n+\dots$,
$G=b_mx^m+\dots$, of degrees $n$ and $m$, is
\begin{equation}
\rest(F,G) = a_n^mb_m^n\prod_{F(\rho)=0}\prod_{G(\eta)=0}(\rho-\eta)
\end{equation}
It is a homogeneous polynomial of degree $m+n$ in the coefficients
of $F$ and $G$, in fact it is homogeneous of degree $m$ in
$a_0,\dots, a_n$  and of degree $n$ in $b_0,\dots, b_m$. Moreover if
$a_i,b_i$ are regarded as having degree $i$ then $\rest(F,G)$ is
homogeneous of degree $mn$. We have
\begin{equation}\label{rest big}
\rest(F,G) = a_n^m\prod_{F(\rho)=0}G(\rho) = (-1)^{mn} b_m^n
\prod_{G(\eta)=0} F(\eta)
\end{equation}

Furthermore, the discriminant of a polynomial $F(x) =
a_nx^n+\dots+a_0$ of degree $n$ may be computed in terms of the
resultant as
\begin{equation}\label{disc in terms of rest}
\disc F = (-1)^{n(n-1)/2}a_{n}^{n-\deg(F')-2} \rest(F,F')
\end{equation}

We apply this to compute the discriminant of $D_f(t) =
\disc(f(x)+t)$, $f(x)=x^n+a_{n-1}x^{n_1}+\dots +a_1x$. The
discriminant $\disc (D_f(t))$ is a polynomial in the coefficients
$a_1,\dots,a_{n-1}$ of $f(x)$.
 We claim that
 the total degree of $\disc D_f(t)$ is $3(n-1)(n-2)$.
 Indeed,  $D_f(t) = \sum_{j=0}^{n-1} b_j t^j$ is a polynomial of
 degree $n-1$ in $t$, and since it is homogeneous of degree $2(n-1)$
 in $t,a_1,\dots, a_{n-1}$ we find that $b_j$ are polynomials of total
 degree  $2(n-1)-j$ in the $a_j$'s. Now $\disc D_f(t) = \sum
 c_r \prod_j b_j^{r_j}$ has total degree $2(n-1)-2=2(n-2)$ in
 the $b_j$'s, that is $\sum r_j = 2(n-2)$,
 and by \eqref{weighted homog}, $\sum_j jr_j = (n-1)(n-2)$. Thus
 the total degree of $\disc D_f(t)$ in $a_1,\dots, a_{n-1}$ is %the maximum, over all $r$, of
\begin{multline*}
\sum_j r_j \deg b_j =  \sum r_j(2(n-1)-j)  = 2(n-1)\sum r_j -\sum
jr_j \\= 2(n-1)\cdot 2(n-2) - (n-1)(n-2) = 3(n-1)(n-2)
\end{multline*}
as claimed.

 Arguing similarly, one sees that the  resultant $\rest(D_f(t),D_{f+\alpha}(t))$  has total degree $3(n-1)^2$
 in the coefficients $a_1,\dots, a_{n-1}$.

 Assume $\gcd(q,n)=1$. Then $f'(t) =
nx^{n-1}+(n-1)a_{n-1}x^{n-2}+\dots$ has degree $n-1$ and by
\eqref{disc in terms of rest}, \eqref{rest big} we find
\begin{equation}\label{expression for D(t) i}
D_f(t)=\disc_x(f(x)+t) =
(-1)^{n(n-1)/2}n^n\prod_{f'(\rho)=0}(t+f(\rho))
\end{equation}
has degree $n-1$, with roots $-f(\rho)$ as $\rho$ runs over the
$n-1$ roots of $f'(x)$.

In the case where $\gcd(q,n)>1$, $f'(t) =-a_{n-1}x^{n-2}+\dots$ has
degree $n-2$ provided that $a_{n-1}\neq 0$, in which case by
\eqref{disc in terms of rest}, \eqref{rest big} we have
\begin{equation}\label{expression for D(t) ii}
D_f(t)=\disc_x(f(x)+t) = (-1)^{n(n-1)/2}a_{n-1}^n
\prod_{f'(\rho)=0}(t+f(\rho))
\end{equation}
which has degree $n-2$ and again has roots $-f(\rho)$ as $\rho$ runs
over the $n-2$ roots of $f'(x)$.

\section{Non-vanishing of the resultant}

\begin{prop}\label{prop:resultant}
Given a nonzero polynomial $\alpha\in \fq[x]$, with $\deg \alpha<n$,
then $a\mapsto \rest(D_f(t),D_{f+\alpha}(t))$ is not the zero
polynomial,
%Let $D_f(t) = \disc_x(f_a(x)+t)$, where $f(x) =f_a(x)=
%x^n+a_{n-1}x^{n-1}+\dots +a_1x$, $\vec a\in \fq^{n-1}$.
%We want to show that the resultant of $D_f(t)$ and $D_{f+\alpha}(t)$ is not
%identically zero,
that is, the polynomial function
\begin{equation}
  R(a):=\rest_t(D_f(t),D_{f+\alpha}(t))\in \Z[\vec a]
\end{equation}
is not identically zero.
\end{prop}
Applying this to $\alpha = \alpha_j - \alpha_i$ for each $j \neq i$
will show that \eqref{Bn small} holds.

\begin{proof}
Write  $ \alpha(x)=A_{n-1}x^{n-1}+\dots +A_0\in \fq[x]$ with $\deg
\alpha<n$.

Let $p$ be the characteristic of $\fq$. Assume first that $p\nmid
n$. Then by \eqref{expression for D(t) i} and \eqref{rest big}, we
find
\begin{equation}
\rest(D_f, D_{f+\alpha}) = n^{2n(n-1)} \prod_{\scriptsize
\begin{tabular}{c}
{$f'(\rho_1)=0$} \\
{$f'(\rho_2)+\alpha'(\rho_2)=0$}
\end{tabular}
} \Big( f(\rho_2)+\alpha(\rho_2)-f(\rho_1) \Big)
\end{equation}

If $p\mid n$ but $a_{n-1} \neq 0$ and $a_{n-1}+A_{n-1}\neq 0$, then
by \eqref{expression for D(t) ii} and \eqref{rest big}, we find
\begin{multline}
\rest(D_f, D_{f+\alpha}) =
a_{n-1}^{n(n-2)}(a_{n-1}+A_{n-1})^{n(n-2)} \\ \times
\prod_{\scriptsize
\begin{tabular}{c}
{$f'(\rho_1)=0$} \\
{$f'(\rho_2)+\alpha'(\rho_2)=0$}
\end{tabular}} \Big(
f(\rho_2)+\alpha(\rho_2)-f(\rho_1) \Big)
\end{multline}
Note that when $a_{n-1} = 0$ or $a_{n-1}+A_{n-1} = 0$, the resultant
$\rest(D_f, D_{f+\alpha})$ is given by different polynomials than in
the above case. However, this might affect at most $2q^{n-2}$
``bad'' $\vec a$'s, which is a negligible amount, and the conclusion
of \eqref{Bn small} remains valid.

In both cases above, the   ``bad'' $\vec a$'s are those for
which there are $\rho_1,\rho_2\in \feq$ such that
\begin{equation}
    f'(\rho_1)=0,\quad f'(\rho_2)=-\alpha'(\rho_2),\quad
f(\rho_2)-f(\rho_1)=-\alpha(\rho_2)
\end{equation}

This is a {\em linear} system of equations for $\vec a\in
\af^{n-1}$, which has the form
\begin{equation}\label{matrix form}
M(\rho)a=b(\rho), \qquad \rho=(\rho_1,\rho_2)
\end{equation}
for a suitable $3\times (n-1)$ matrix $M(\rho)$ and vector
$b(\rho)\in \af^3$. Thus over $\feq$, the solutions of $R(\vec a)=0$
are precisely those $\vec a\in \feq^{n-1}$  which satisfy the system
\eqref{matrix form} for some $\rho\in \feq^2$.

We consider the affine variety (possibly reducible) defined by these
equations
\begin{equation}
Z=\{(\rho,a)\in \af^2\times \af^{n-1}: M(\rho)a=b(\rho) \} \subset
\af^2\times \af^{n-1}
\end{equation}
Let $\phi:Z\to \af^{n-1}$ be the restriction to $Z$ of the
projection $\af^2\times \af^{n-1}\to \af^{n-1}$, and $\pi:Z\to
\af^2$ the restriction to $Z$ of the projection $\af^2\times
\af^{n-1}\to \af^{2}$.
\begin{equation}
\xymatrix{ &{Z\subset \af^2\times \af^{n-1}} \ar[dl]_{\pi}
\ar[dr]^{\phi}
&\\
{\af^2}  && {\af^{n-1}} }
\end{equation}
From the above, the solution set of $R(\vec a)=0$ is precisely
$\phi(Z)$.

We will show that $Z$ has dimension $n-2$, and hence the dimension
of $\{R=0\}=\phi(Z)$ cannot exceed $n-2$ and hence is not all of
$\af^{n-1}$. Thus $R$ is not the zero polynomial, proving
Proposition~\ref{prop:resultant}.

To do so, we study the dimensions of the fibers $\pi^{-1}(\rho)$,
which are affine linear subspaces. We first assume that $n>3$. In
this case we will show that $\pi(Z)$ is dense in $\af^2$ and
generically, that is if $\rho_1\neq \rho_2$, the fibers
$\pi^{-1}(\rho)$ have dimension $n-4$. Moreover there are at most
$\deg \alpha$ non-generic fibers, each of dimension $n-2$.  This
will show that $\dim Z=n-2$.

We rewrite the system \eqref{matrix form} as:
\begin{equation}\label{expanded matrix form}
\begin{split}
 %&a_{n-1}(n-1)\rho^{n-2} +
\dots  +3a_3\rho_1^2+2a_2\rho_1+a_1&=-n\rho_1^{n-1}  \\
%&a_{n-1}(n-1)\rho_2^{n-2} +
\dots +3a_3\rho_2^2+2a_2\rho_2+a_1
&=-\alpha'(\rho_2)-n\rho_2^{n-1} \\
%&a_{n-1}\srho_2^{n-1}+
\dots +a_3
(\rho_2^3-\rho_1^3)+a_2(\rho_2^2-\rho_1^2)+a_1(\rho_2-\rho_1)&=
-\alpha(\rho_2)-(\rho_2^n-\rho_1^n)
\end{split}
\end{equation}
  To find the rank of the matrix $M(\rho)$ we compute that
\begin{equation}
\det \begin{pmatrix}
  3\rho_1^2&2\rho_1&1  \\
3\rho_2^2&2\rho_2&1\\
\rho_2^3-\rho_1^3& \rho_2^2-\rho_1^2&\rho_2-\rho_1
\end{pmatrix} = (\rho_1-\rho_2)^4
\end{equation}
and thus $M(\rho)$ has full rank $3$ unless $\rho_1=\rho_2$, and so
the generic fibers $\pi^{-1}(\rho)$ have dimension $n-1-3=n-4$.

In the non-generic case $\rho_1=\rho_2$, the matrix has rank one and
we need $\alpha'(\rho_2)=0=\alpha(\rho_2)$, which constrains us to
have at most finitely many fibers (the number bounded by $\deg
\alpha$/2), each of which has dimension $n-1-1=n-2$.

Finally, the cases $n=2,3$ are handled similarly, except that  the
image of the map $\pi:Z\to \af^2$ is no longer dense, due to
algebraic conditions constraining $\rho_1,\rho_2$. We omit the
(tedious) details.
\end{proof}

\section{Non-vanishing of the discriminant}

% To show that the condition for being in $A_n$ is not always
% satisfied,  it suffices to show that outside a set of $\vec a$'s
% of codimension at least one in the parameter space $\af^{n-1}$,
% $D_f(t)$ is of positive degree, and is square-free, that is with nonzero discriminant.
% Note that $\disc(D_{f})$ is of total degree $\leq 3(n-1)(n-2)$ in
% $\vec a$.

We wish to show that the condition for being in $A_n$ is not always
satisfied. Without loss of generality, we can  assume $\alpha_i = 0$. We
first study a couple of small degree cases.

For $n=2$, $\disc(x^2+ax+t) = a^2-4t$ is linear and hence has no
repeated roots (recall $q$ is odd), hence $A_n$ is empty. When
$n=3$, we have
\begin{equation}
D_f(t)=\disc_x (x^3+ax^2+bx+t) =\left(a^2 b^2-4 b^3\right)+ \left(18
a b-4 a^3\right)t-27 t^2
\end{equation}
If $3\mid q$ then $D_f(t) =\left(a^2 b^2-4 b^3\right)-4a^3t$ has
degree one for $a\neq 0$; if $3\nmid q$ then $D(t)$ has degree $2$
and we compute that
\begin{equation}
\disc_t \disc_x (x^3+ax^2+bx+t) = -16 (a^2-3b)^3
\end{equation}
which is clearly not identically zero.  So we may assume
$n\geq 4$. %In the case that $\gcd(q,n(n-1))=1$,

% To handle the case when $\gcd(q,n(n-1))>1$, we need to argue
% differently.

%Thus \eqref{An small} holds.

\subsection{}
Similarly to our approach in the previous section, it suffices to
show that outside a set of $\vec a$'s of codimension at least one in
the parameter space $\af^{n-1}$, $D_f(t)$ is of positive degree, and
is square-free, that is with nonzero discriminant.

We conclude from \eqref{expression for D(t) i}, \eqref{expression
for D(t) ii} that if $n\geq 4$ and $\vec a$ is in the ``bad'' set
(but $a_{n-1}\neq 0$ if $\gcd(n,q) \neq 1$), then at least one of
the following occurs:
\begin{itemize}
\item
There is some $\rho\in \feq$  for which $f'(x)$ has a double zero at
$x=\rho$, that is, there is some $\rho \in \feq$ for which
\begin{equation}\label{first system}
 f'(\rho)=0,\quad f''(\rho)=0
\end{equation}
\item There are two distinct $\rho_1\neq \rho_2$
so that $f(\rho_1)=f(\rho_2)$ and  so that $f'(x)$ vanishes at both
$x=\rho_1$ and $x=\rho_2$, that is
\begin{equation}\label{second system}
 f'(\rho_1)=0,\quad
  f'(\rho_2)=0 ,\quad  f(\rho_1)=f(\rho_2)
\end{equation}
\end{itemize}
We want to show that the set of $\vec a\in \feq^{n-1}$ which solve
at least one of \eqref{first system}, \eqref{second system} has
dimension at most $n-2$.
\subsection{}
We first look at $f$ for which \eqref{first system} happens. This
gives a pair of equations for $\vec a\in \feq^{n-1}$:
\begin{equation}\label{rho system}
\begin{split}
\dots + 2\rho a_2+a_1 &= -n\rho^{n-1}\\
\dots+ 2a_2+0 &=-n(n-1)\rho^{n-2}
\end{split}
\end{equation}
Defining
\begin{equation}
W=\{(\rho,\vec a)\in \af^1\times \af^{n-1}: \eqref{first system}
\mbox{
 holds}\}
\end{equation}
we have a fibration of $W$ over the $\rho$ line $\af^1$ and a map
$\phi: W\to \af^{n-1}$, the restriction of the projection
$\af^1\times \af^{n-1}\to \af^{n-1}$,
\begin{equation}
\xymatrix{ &{W\subset \af^1\times \af^{n-1}} \ar[dl]_{\pi}
\ar[dr]^{\phi}
&\\
{\af^1}  && {\af^{n-1}} }
\end{equation}
and the solutions of \eqref{first system} are precisely $\phi(W)$.

The system \eqref{rho system} is nonsingular (rank $2$) and hence
$\pi:W\to \af^1$ is surjective and for each $\rho$ the dimension of
the solution set is $n-1-2=n-3$.
 we find that $\dim W=n-2$ and hence $\dim \phi(W)\leq n-2$.

\subsection{}
Next we consider the system \eqref{second system} which given
$\rho_1\neq \rho_2$ is a linear system for $\vec a\in \feq^{n-1}$,
of the form
\begin{equation}\label{Third system}
\begin{split}
 \dots +3\rho_1^2a_3+2\rho_1 a_2+a_1 &= -n\rho_1^{n-1}\\
 \dots +3\rho_2^2a_3+2\rho_2 a_2+a_1 &= -n\rho_2^{n-1}\\
\dots +(\rho_2^3-\rho_1^3)a_3+(\rho_2^2-\rho_1^2)a_2 +
(\rho_2-\rho_1)a_1 &= -\rho_2^n+\rho_1^n
\end{split}
\end{equation}

This system shares the matrix part of \eqref{expanded matrix form},
and hence has rank 3 for every $\rho_1 \neq \rho_2$. Thus the
arguments of the previous section show that
\begin{equation}
\{\vec a \in \af^{n-1} : \exists \rho_1 \neq \rho_2 \mbox{\; s.t.
\eqref{second system}  holds}\}
\end{equation}
is of dimension at most $n-2$. This shows that \eqref{An small}
holds, thus concluding the proof of Proposition~\ref{Prop Gcomp}.

\end{document}